\documentclass{mpi2015-cscpreprint}

\usepackage[american]{babel}
\usepackage{amssymb}
\usepackage{amsthm,amsmath}
\usepackage{algorithm}%
\usepackage{algorithmicx}%
\usepackage{algpseudocode}%
\usepackage{graphicx}
\usepackage{appendix}
\usepackage{multirow}

\newtheorem{theorem}{Theorem}
\newtheorem{remark}{Remark}
\newtheorem{assumption}{Assumption}

\numberwithin{equation}{section}

\begin{document}

\title{A Model-Free Extremum Seeking Controller with Application to Tracking a Nonlinear Chemical Reaction}

\author[1,3]{Alexander~Zuyev}
\author[2,3]{Victoria Grushkovska}
\affil[1]{Max Planck Institute for Dynamics of Complex Technical Systems, 39106 Magdeburg, Germany}
\affil[2]{Department of Mathematics, University of Klagenfurt, 9020 Klagenfurt am W\"orthersee, Austria}
\affil[3]{Institute of Applied Mathematics and Mechanics, National Academy of Sciences of Ukraine}

\keywords{extremum seeking, Lie bracket, nonlinear chemical reaction}

\msc{93C40, 49N30, 93D21, 93C95}

\abstract{In this paper, we develop the extremum-seeking approach to generate admissible trajectories in a neighborhood of a given reference curve in the state space. The cost function of the problem represents the distance between the current system state and the reference curve, which is parameterized as a function of time. Such reference curves naturally arise as optimal trajectories in isoperimetric optimization problems for nonlinear chemical reactions, where the objective is to maximize the average reaction product over a given period. We apply the proposed extremum seeking control design to a nonisothermal reaction model and illustrate the resulting tracking errors through numerical simulations.}

\maketitle

\section{Introduction}
Extremum seeking is known to be an efficient control methodology for various dynamic optimization problems, requiring only mild assumptions about the governing equations and the structure of the cost function (a recent survey in this area is available in~\cite{scheinker2024100}).
Model-free extremum seeking algorithms are particularly important in control design for complex systems with uncertain dynamics, see, e.g.,~\cite{Kr00,Har13,Guay15}.
Such algorithms can be especially beneficial for optimizing chemical and biochemical processes with partially defined kinetics or uncertain kinetic parameters~\cite{dochain2011extremum,lopez2024line}.
In addition to the constrained steady-state optimization problem considered in~\cite{dochain2011extremum},
it has been shown that forced periodic operations can improve the performance of chemical reactions compared to a steady-state regime (see, e.g.,~\cite{silveston2012periodic} and references therein).
A rigorous analysis of the related isoperimetric optimal control problem with periodic boundary conditions for a continuous stirred-tank reactor (CSTR) model has been conducted in~\cite{CES2017} using the Pontryagin maximum principle.
Additionally, periodic trajectories for a class of bang-bang controls have been constructed in~\cite{AMM2019} using the Chen--Fliess series expansion.
 These results have been illustrated with the hydrolysis reaction of acetic anhydride, and an experimental demonstration of the product yield for this reaction is presented in~\cite{felischak} for the case of sinusoidal input modulations.
In analyzing the experimental data, the authors of~\cite{felischak} noted that
the reactor does not appear to be fully adiabatic, leading to deviations of measured signals from their predicted values and causing a loss in mean product yield.
This outcome motivates the development of adaptive control techniques to refine the design of required periodic trajectories, particularly by accounting for uncertain aspects of the reactions under consideration.

Motivated by the above conclusion, we focus here on extending the extremum-seeking methodology with a model-free component to generate controls based solely on output measurements, ensuring that deviations from a specified reference trajectory remain practically small.

The main contribution of this paper is twofold. On one hand, we  extend the extremum seeking control approach for dynamic systems to the case of a time-varying cost function.
On the other hand, this methodology is applied to tracking a nonlinear chemical reaction model in the proximity of a given curve in the state space, based on a scalar cost that measures only the magnitude of the tracking error.

\section{Extremum seeking problem}
Consider the system
\begin{equation}\label{sys_dyn}
\begin{aligned}
  &\dot x=f(x,u),\\
  &y=   h(t,x),
\end{aligned}
\end{equation}
where  $x\in D_x\subseteq\mathbb R^{n_x}$ is the state vector, $x(t_0)=x^0\in D_x$ ($t_0\in\mathbb R^+$), $u\in D_u\subseteq \mathbb R^{n_u}$ is the control vector,  $y\in D_y\subseteq \mathbb R^+$ represents the output of the system, $f\in C^1(D_x\times D_u;\mathbb R^{n_x})$, $h \in C^1(\mathbb R^+\times D_x; D_y)$, $n_x,n_u,n_y\in\mathbb N=\{1,2,...\}$. %

In this paper, we consider the problem of finding a control law $u$ that optimizes the steady-state of the system in the practical sense that the system's output can be made smaller than a given tolerance level $\rho>0$, i.e.
$$
h(t,x)\le \rho\text{ for large enough }t.
$$
 Such problems arise, for example,
 in trajectory tracking problems,
 where  $y=\|x(t)-x^*(t)\|^2$ represents the tracking error, or
  in output regulation problems. Unlike many existing results in this area, we assume that only the values of  $h(t,x)$ are available for control design. Furthermore, we assume that the output may not identically vanish  at a constant point  $x^*$,  but rather along a curve  $x^*(t)$.

To address this problem, we adopt extremum-seeking methods based on the Lie bracket framework~\cite{Durr17,SK17,GZE18,Gu18}.
As it is shown in~\cite{Durr17,GE21}, in case of time-invariant output $y=h(x)$, the dynamic controller
$$
\dot u=\frac{\gamma}{\eta\sqrt\varepsilon}\sum_{j=1}^{n_u}\sqrt j\Big(h(x)\cos\Big(\frac{jt}{\eta\varepsilon}\Big)+\sin\Big(\frac{jt}{\eta\varepsilon}\Big)\Big)e_j
$$
can be used to ensure the singular practical uniform asymptotic stability of a constant optimal state for system~\eqref{sys_dyn}. Here,  $e_j$ denotes the $j$-th unit vector in $\mathbb R^{n_u}$, and $\varepsilon,\eta$ are some positive parameters. By \emph{singular practical uniform asymptotic stability} we mean that $\varepsilon$ should be sufficiently small  and $\eta$ sufficiently large  to ensure the desired proximity of the solutions to the optimal state. For precise definitions, see~\cite{Durr17,GE21}.
 Simply saying,  the role of the parameter $\varepsilon$ is to ensure that the trajectories $\tilde u(t)$ of  the reduced system
 $$
\dot {\tilde u}_j=\frac{1}{\eta\sqrt\varepsilon}\sum_{j=1}^{n_u}\sqrt j\Big(h(\ell(\tilde u))\cos\Big(\frac{jt}{\eta\varepsilon}\Big)+\sin\Big(\frac{jt}{\eta\varepsilon}\Big)\Big)e_j
$$
 tend asymptotically to a neighborhood of
 $$
 u^*={\textrm arg} \min\limits_{u\in D_u} h(\ell(u)),
 $$
 where $\ell(u)$ is the steady-state map~\cite{Durr17,GE21} of system~\eqref{sys_dyn}.
 The parameter
 $\eta$  ensures that the extremum seeking system evolves at a significantly slower rate compared to the controlled plant.
In this paper, we extend the applicability of this approach to the case of time-varying  $h(t,x)$.

\section{Main results}

\subsection{Control design}

In order to apply the approach described in the previous section, we need to formulate certain assumptions. First of all, suppose that there exists a so-called steady-state map $\ell:D_u\to D_x$  such that, for each constant $\bar u\in D_u$,
and for any $x^0\in D_x$, $t_0\in\mathbb R^+$, the solutions $x(t)$ of system~\eqref{sys_dyn} with $u(t)\equiv\bar u$ and $x(t_0)=x^0$ satisfy
$$
x(t)\to \ell(\bar u)\text{ as }t\to+\infty.
$$
More precisely, we require the following property.
\begin{assumption}\label{ass_ell}
A function  $f\in C^1(\mathbb R^{n_x}\times \mathbb R^{n_u};\mathbb R^{n_x})$, and there exists a function $\ell\in C^2(\mathbb R^{n_u}; \mathbb R^{n_x})$  such that, for each fixed $u\in D_u$, $f(x,u)=0$ if and only if $x=\ell(u)$.
\end{assumption}
Moreover, it is important to have the asymptotic stability of the state $\ell(u)$ of system~\eqref{sys_dyn} for each fixed (``frozen'') $u\in D_u$. This is ensured, in particular, if system~\eqref{sys_dyn} admits a Lyapunov function for each fixed $u$, as stated in the next assumption. However, further relaxation of this requirement are possible, as it is shown in~\cite{GE21,GZ20_IFAC}.
\begin{assumption}\label{ass_V}
There exist  a function $V\in C^1(D_u\times D_x;\mathbb R)$ and positive constants $\sigma_{11},\sigma_{12},\sigma_{21},\sigma_{22}$ such that, for all $x\in
      D_x$, $u\in D_u$,
  $$
  \begin{aligned}
    \sigma_{11}^2\|x-l(u)\|^2&\le {V(x,u)}\le \sigma_{12}^2\|x-l(u)\|^2\\
    \frac{\partial V(x,u)}{\partial x}f(x,u)&\le-\sigma_{21}{V(x,u)},\\
     \Big\|\frac{\partial V(x,u)}{\partial u}\Big\|^2&\le\sigma_{22}^2{V(x,u)}.
  \end{aligned}
  $$
\end{assumption}
The next assumptions concern properties of the output map.
\begin{assumption}\label{ass_h}
The functions $h(t,x)$ and $u^*(t)$ satisfy the following properties with some $\alpha_{11},\alpha_{12}$, $\alpha_{21},\alpha_{22},\alpha_3$, and $L_h>0$ for all $x,\tilde x\in D_x$, $u\in D_u$, $t\ge0$:
\begin{enumerate}
    \item there exists an $L_h>0$ such that $$|\sqrt{h(t,x)}-\sqrt{h(t,\tilde x)}|\le L_h\|x-\tilde x\|$$  for all $x,\tilde x\in D_x$, $t\ge0$;
    \item there exist $\alpha_{11},\alpha_{12},\alpha_{21},\alpha_{22},\alpha_3>0$ such that
    $$
      \begin{aligned}
\alpha_{11}\|u-u^*(t)\|\le \sqrt{ h(t,\ell(u))} &\le\alpha_{12}\|u-u^*(t)\|,\\
\alpha_{21}\sqrt{ h(t,\ell(u))}\le \|\nabla_u  h(t,\ell(u))\|&\le\alpha_{22} \sqrt{ h(t,\ell(u))},\\
\Big\|\frac{\partial^2 h(t,\ell(u))}{\partial u^2}\Big\| &\le\alpha_3
  \end{aligned}
    $$
    for all $u\in D_u$, $t\ge0$;
    \item there exists a  $\nu\ge0$ such that
    $$
\sup\limits_{t_1,t_2\in\mathbb R^+} \|u^*(t_1)-u^*(t_2)\|\le  \nu.
    $$
\end{enumerate}
\end{assumption}
To simplify the presentation, for the remainder of this paper, we define
$$
D_x=\bigcup\limits_{t\ge0}B_{\Delta_x}(\ell(u^*(t))),\; D_u=\bigcup\limits_{t\ge0}B_{\Delta_u}(u^*(t))
$$
for some $\Delta_x,\Delta_u\in(0,+\infty]$.

We exploit the extremum seeking approach. Accordingly, we construct the following dynamical extension for system~\eqref{sys_dyn}:
\begin{equation}\label{cont_dyn}
\dot u_j=\frac{2\gamma}{\eta\sqrt\varepsilon}\sqrt{\pi  jh(t,x)}\sin\left(\ln h(t,x)+\frac{2\pi j t}{\eta\varepsilon}\right),
\end{equation}
with $j=1,\dots,n_u$,  $u(t_0)=u^0\in D_u$, $\varepsilon,\eta,\gamma>0$.
As the functions $$g_1(t,u)=\sqrt{h(t,\ell(u))}\sin\left(\ln h(t,\ell(u))\right)\;\; \text{and} \;\; g_2(t,u)=\sqrt{h(t,\ell(u))}\cos\left(\ln h(t,\ell(u))\right)$$ are not differentiable if $h(t,\ell(u))=0$, we  exploit the following regularity assumption from~\cite{GZE18,GE21,GMZME18}:
\begin{assumption}
The functions $g_1(t,u)$ and $g_2(t,u)$ are twice continuously differentiable for all $t\in\mathbb R^+$, $u\in D_u\setminus \bigcup\limits_{t\ge0}u^*(t)$. Furthermore, the functions $\mathcal L_{g_{j_2}}g_{j_1}(t,u)$, $\mathcal L_{g_{j_3}}\mathcal L_{g_{j_2}}g_{j_1}(t,u)$ are continuous, and the first-order partial derivatives of $g_{j_1}(t,u)$, $\mathcal L_{g_{j_2}}g_{j_1}(t,u)$ with respect to $t$ are bounded  for all $t\in\mathbb R^+$, $u\in D_u\setminus \bigcup\limits_{t\ge0}u^*(t)$, $j_1,j_2,j_2\in\{1,2\}$. Here, $$\mathcal L_{g_{j_2}}g_{j_1}(t,u):=\lim\limits_{s\to 0}\dfrac{1}{s}\big(g_{j_1}(t,u+sg_{j_2}(t,u))-g_{j_1}(t,u)\big).$$
\end{assumption}

\subsection{Stability properties}
\begin{theorem}
Given system~\eqref{sys_dyn}--\eqref{cont_dyn}, let  Assumptions 1--4 hold. Then, for any $\delta_x,\delta_u>0$ and $\rho>\nu$, there exist $\bar\varepsilon>0$ and $\bar\eta(\varepsilon)>0$ such that, for any $\varepsilon\in({\nu}/{\gamma\rho},\min\{\bar\varepsilon,{1}/{\gamma}\})$, $\eta\in(\eta(\varepsilon),\infty)$, the solutions of system~\eqref{sys_dyn}--\eqref{cont_dyn} with initial conditions $(x^0,u^0)\in\overline{B_{\delta_{x}}(\ell(u^*(t)))\times B_{\delta_u}(u^*(t))}$ satisfy the property
$$
\|x(t)-\ell(u^*(t))\|+\|u(t)-u^*(t)\|\le \rho\text{ for all }t\in[t_f,+\infty),
$$
for some $t_f\ge 0$.
\end{theorem}

The proof of Theorem~1 is presented in the Appendix.
\begin{remark}
Unlike conventional extremum seeking results, which assume control frequencies to be infinitesimally small, we require $\varepsilon$ to be bounded from below by a certain value. This requirement arises from the following consideration.   On the one hand, subsystem~\eqref{cont_dyn} must evolve sufficiently slowly relative to subsystem~\eqref{sys_dyn}, as this separation of time-scales is a fundamental characteristic of extremum seeking for dynamical systems~\cite{Durr17,GE21}. This condition is ensured by choosing a sufficiently large $\eta$. On the other hand, $\varepsilon$ must not be too small to guarantee sufficiently close tracking of the curve  $\ell(u^*(t))$. To achieve this, we impose specific constraints on  $\varepsilon $ and $\gamma$. Furthermore, the parameter $\nu$ in Assumption~3 must be sufficiently small to ensure close convergence to the target curve. In future work, we plan to revisit these restrictions and explore alternative attractivity conditions.
\end{remark}

 \section{Application to a controlled chemical reaction}
In this section, we apply the proposed control scenario to a non-isothermal chemical reaction model of the type ``$A \to$ Product'' considered in~\cite{CES2017}.
The model is represented by:
\begin{equation}\label{cstr}
\dot x = f(x,u),\;\; x=\begin{pmatrix}x_1 \\ x_2 \end{pmatrix}\in D\subset {\mathbb R}^2,\; u=\begin{pmatrix}u_1 \\ u_2 \end{pmatrix} \in U\subset {\mathbb R}^2,
\end{equation}
where
$$
f(x,u)=\begin{pmatrix}
- \phi_1 x_1 + k_1 e^{-\varkappa} - (x_1+1)^{\bar n} e^{-\varkappa/(x_2+1)} + u_1\\
- \phi_2 x_2 + k_2 e^{-\varkappa} - (x_1+1)^{\bar n} e^{-\varkappa/(x_2+1)} + u_2
\end{pmatrix},
$$
$$
D=\{x\in{\mathbb R}^2\,|\,x_1>-1,\;x_2>-1\},
$$
$$
U=\{u\in{\mathbb R}^2\,|\, u_1\in [u_1^{min},u_1^{max}],\; u_2\in [u_2^{min},u_2^{max}]\}.
$$
Here,
$x_1$ and $x_2$ represent the concentration of $A$ and the temperature in the reactor, respectively, while
 $u_1$ and $u_2$ control the inlet concentration and inlet temperature.
We assume that the state variables and controls of this system are dimensionless deviations from a steady-state operating mode of a continuous steered-tank reactor (CSTR).
This physical steady-state corresponds to the trivial equilibrium $x_1=x_2=0$ with $u_1=u_2=0$ for system~\eqref{cstr}.
The reaction order $\bar n$ and parameters $\phi_1$, $\phi_2$, $k_1$, $k_2$, $\varkappa$ are chosen as in~\cite{CES2017}:
\begin{equation}\label{parameters}
\begin{aligned}
&{\bar n}=1,\;\phi_1=\phi_2=1,\\
&k_1=5.819\cdot 10^7,k_2=-8.99\cdot 10^5,\; \varkappa=17.77, \\ &u_1^{max}=-u_1^{min} = 1.798,\; u_2^{max}=-u_2^{min}=0.06663.
\end{aligned}
\end{equation}
The above parameter values correspond to the hydrolysis reaction
$$\mathrm{(CH_3CO)_2 O + H_2 O \to 2\, CH_3 COOH}
$$ with excess of water.

The Jacobian matrix of $f(x,u)$ at the equilibrium point $(x,u)=(0,0)$ is given by:
\begin{equation}\label{A_chem}
\begin{aligned}
\tilde A &= \left.\frac{\partial f(x,u)}{\partial x}\right|_{(x,u)=(0,0)} \\
&=\begin{pmatrix}-\phi_1-k_1 e^{-\varkappa} &-k_1 \varkappa e^{-\varkappa}\\
-k_2 e^{-\varkappa} & -\phi_2-k_2\varkappa e^{-\varkappa}
\end{pmatrix}.
\end{aligned}
\end{equation}
Using the parameter values from~\eqref{parameters}, we obtain:
$$
\tilde A \approx \begin{pmatrix}-2.115412260 & -19.82087587\\
0.01723243894 & -0.6937795600
\end{pmatrix}.
$$
The eigenvalues of $\tilde A$ are $\lambda_1\approx -1$ and $\lambda_2\approx -1.809$. This implies that the matrix $\tilde A$ is nonsingular.
Therefore, by the implicit function theorem, for each $u$ in some neighborhood $\cal U$ of $0$ in ${\mathbb R}^2$, there exists a unique $x=\ell (u)\in D$ such that
$f(\ell (u),u)=0$. Since $f$ is of class $C^1$, the Jacobian matrices
$$
 \left.\frac{\partial f(x,u)}{\partial x}\right|_{(x,u)=(\ell (u),u)}
$$
are Hurwitz for each $u$ from some neighborhood of zero ${\cal U}_0\subseteq \cal U$.

Control system~\eqref{cstr} admits a periodic solution $x^*(t)$ with the period $T=100$ for the following control functions:
\begin{equation}\label{u_trig}
u^*_j(t)=u_j^{min}\sin\left(\dfrac{2\pi t}{T}\right),\; j=1,2.
\end{equation}
The initial values of this solution are: $x^*_1(0)\approx -0.065$ and $x^*_2(0)\approx 0.008$.
Trigonometric periodic controls have been experimentally tested on the hydrolysis reaction of the type described by system~\eqref{cstr} in the paper~\cite{felischak}, where it was shown that they lead to a better product yield compared to the steady-state operation.
Taking into account that the mathematical model of this reaction may have uncertainties, we consider here the problem of tracking the above-introduced
reference curve $\{x^*(t)\}$ using only the information on current deviations in the concentration and temperature:
\begin{equation}\label{output_chhem}
y=h(t,x) = \|x-x^*(t)\|^2.
\end{equation}
For numerical simulations, we apply controller~\eqref{cont_dyn} with the following parameters:
$$
 \gamma=150,\,\varepsilon=0.001,\,\eta=1.
$$
It should be emphasized that controller~\eqref{cont_dyn} is model-free, and its implementation requires only the measurement of the output signal~\eqref{output_chhem} without requiring knowledge of the exact reaction model.
Simulation results for system~\eqref{cstr} with output~\eqref{output_chhem} and controller~\eqref{cont_dyn} are shown in Fig.~1, a).

In addition to the trigonometric inputs, we also test the proposed extremum seeking strategy for the following bang-bang controls:
\begin{equation}\label{u_bang}
u^*_j(t)=u_j^{min}{\textrm sign}\,\sin\left(\dfrac{2\pi t}{T}\right),\; j=1,2,
\end{equation}
and the corresponding periodic solution $x^*(t)$ of system~\eqref{cstr}.
The results of numerical simulations for this trajectory are reflected in Fig.~1, b).

 \begin{figure*}[thpt]\label{ex_co}
   \begin{minipage}{1\linewidth}
  \centering
  \includegraphics[width=0.49\linewidth]{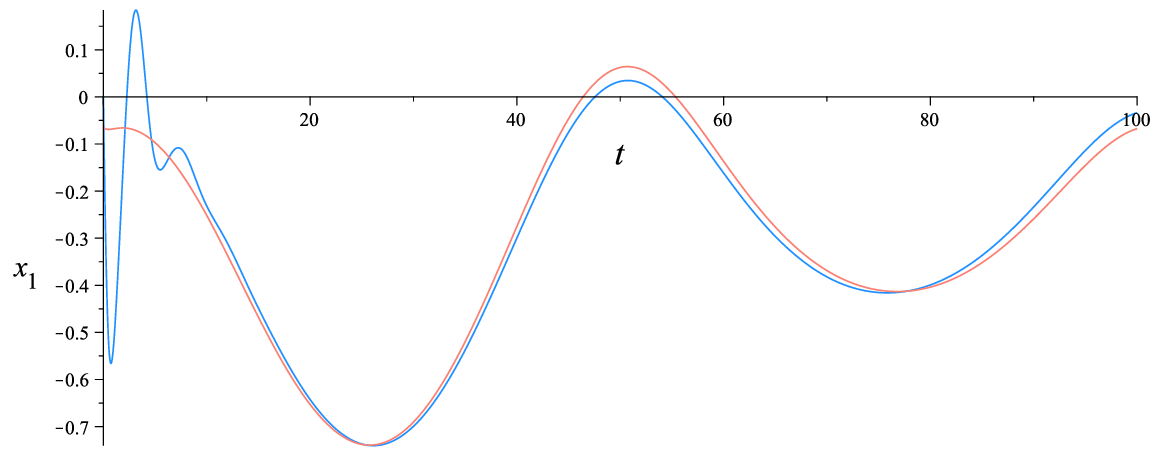}\hfill
\includegraphics[width=0.49\linewidth]{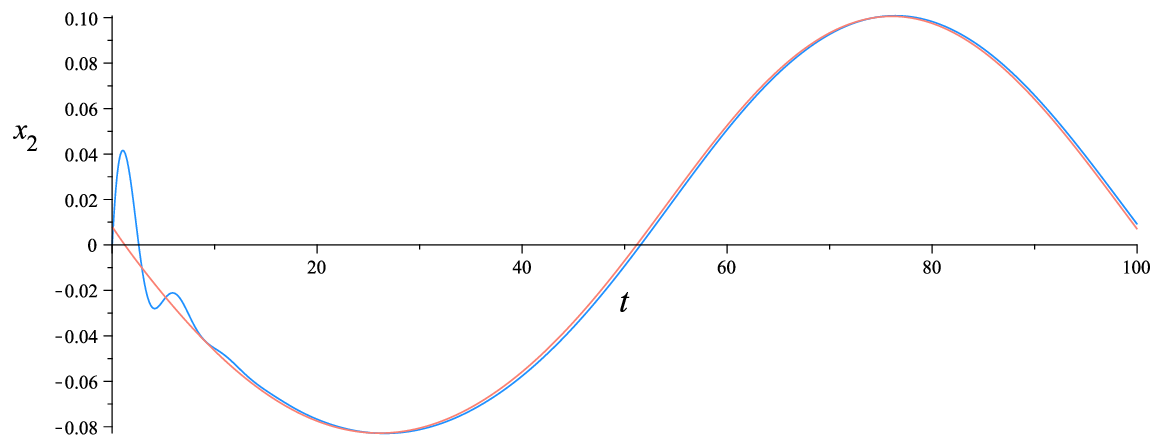}\\
a)\\
\includegraphics[width=0.49\linewidth]{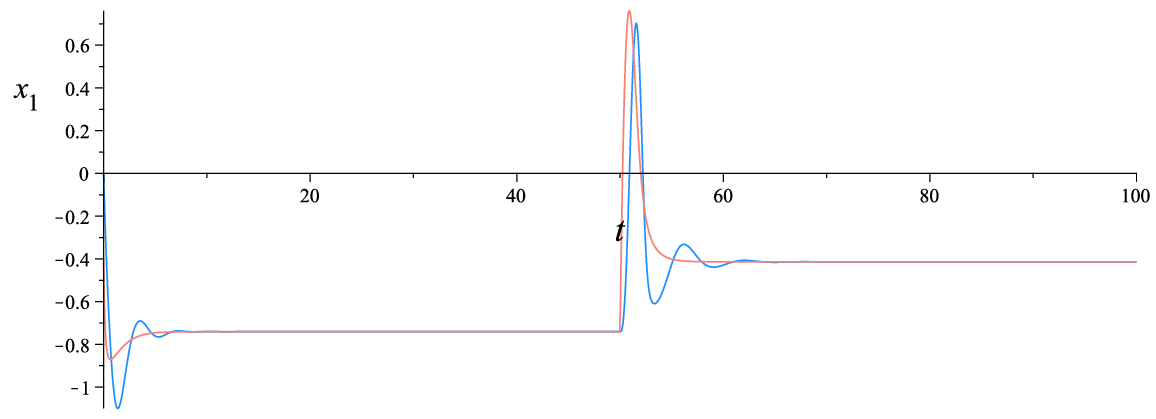}\hfill
\includegraphics[width=0.49\linewidth]{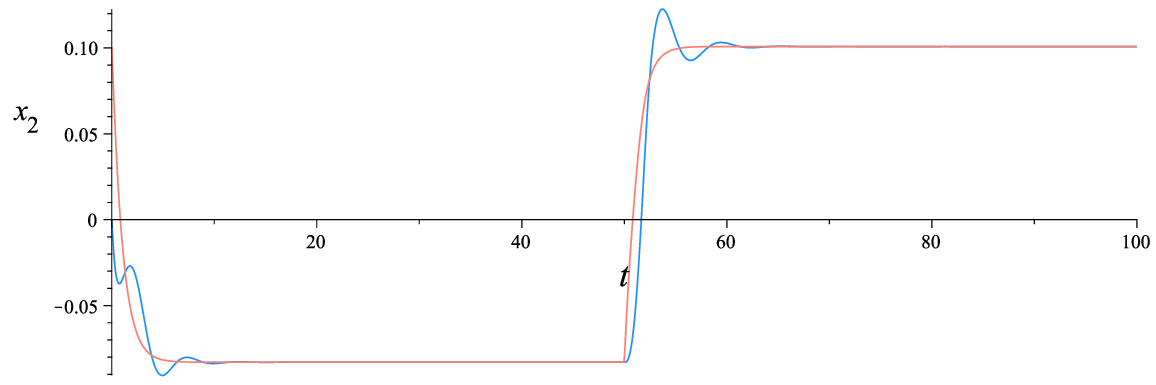}\\
b)
 \end{minipage}\hfill
\caption{Time plots of the $x_1$ (left) and $x_2$ (right) components of system~\eqref{cstr} under  controller~\eqref{cont_dyn}, tracking a reference trajectory generated by the control inputs~\eqref{u_trig} in a) and~\eqref{u_bang} in b). In all plots, the reference trajectory is shown in red, and the system response is shown in blue. }
 \end{figure*}

 \section{Conclusion}
The obtained results justify the feasibility of tracking a given curve using the two-scale extremum seeking controller with a model-free component, based solely on measurements of the tracking error as a function of time. The proof of practical convergence in Theorem~1 requires the regularity of the reference curve and the corresponding input signal. However, the numerical example demonstrates that the proposed controller can also ensure the practical convergence of solutions for a realistic chemical reaction model with discontinuous inputs.

A tradeoff between the convergence rate of the reduced system and the approximation of the original controlled dynamics can be achieved through a suitable choice of the time-scale parameters $\varepsilon$ and $\eta$,  where the theoretical bounds for these parameters are established in the proof of Theorem~1.

\newpage
\appendix
\section{Appendix: Proof of Theorem~1}

The proof  extends the approach proposed in~\cite{GE21} to the tracking problem.

Similarly to the results of~\cite{GE21}, it can be shown that, for any $\delta_x\in(0,\Delta_x)$, $\delta_y\in(0,\Delta_y)$, there exist an  $\varepsilon_1>0$ and an $\eta_0:(0,\infty)\to(0,\infty)$ such that, for any $\varepsilon\in(0,\varepsilon_1)$, $\eta\in(\eta_0(\varepsilon),\infty)$, and any  $(x^0,u^0)\in\overline{B_{\delta_{x}}(\ell(u^*(t)))\times B_{\delta_u}(u^*(t))}$, the solutions $(x(t),u(t))$ of system~\eqref{sys_dyn}--\eqref{cont_dyn} with initial condition $x(0)=x^0$, $u(0)=u^0$ are well-defined in $D_x\times D_u$ for all $t\ge0$. We omit this step because of the page limits and assume in the sequel that $\varepsilon\in(0,\varepsilon_0)$, $\eta\in(\eta_0(\varepsilon),\infty)$.

\vskip1em

\emph{\textbf{Step 1.} Estimating the maximal deviations between the $x$ and $u$ components of the solution of system~\eqref{sys_dyn}--\eqref{cont_dyn} with $(x^0,u^0)\in\overline{B_{\delta_{x}}(\ell(u^*(t)))\times B_{\delta_u}(u^*(t))}$.}

Using Assumptions~2--3 and H\"older's inequality, we get the following estimates:
$$
\begin{aligned}
 &\sqrt{h(t,\ell(u(t)))}\le \alpha_{12}\|u(t)-u^*(t)\|\\
& \le \alpha_{12}\left(\|u^0-u^*(0)\|+\|u^*(t)-u^*(0)\|+\|u(t)-u^0\|\right)\\
 & \le \alpha_{12}\left(\|u^0-u^*(0)\|+\nu+\int\limits_0^t\|\dot u(s)\|ds\right)\\
& \le \alpha_{12}\left(\|u^0-u^*(0)\|+\nu+\frac{\gamma c_u}{\eta\sqrt\varepsilon}\int\limits_0^t \sqrt{h(s,x(s))}ds\right),
\end{aligned}
$$
where $c_u=\gamma \sqrt{2\pi n_u(1+n_u)}$.
Thus,
$$
\begin{aligned}
&\sqrt{h(t,x(t))}\le \sqrt{h(t,\ell(u(t)))}+L_h\|x(t)-\ell(u(t))\|\\
&\le \alpha_{12}\left(\|u^0-u^*(0)\|+\nu+\frac{c_u}{\eta\sqrt\varepsilon}\int\limits_0^t \sqrt{h(s,x(s))}ds\right)\\
&+\frac{L_h}{\sigma_{11}}\sqrt{V(x(t),u(t))}.
\end{aligned}
$$
From the Gr\"onwall--Bellman inequality, this implies
\begin{equation}\label{est_h}
\begin{aligned}
&\sqrt{h(t,x(t))}\le \alpha\left(\|u^0-u^*(0)\|+\nu\right)\\
&+\frac{L_h}{\sigma_{11}}\sqrt{V(x(t),u(t))}+\frac{\alpha L_hc_u}{\sigma_{11}\eta\sqrt\varepsilon}\int\limits_0^t \sqrt{V(x(s),u(s))}ds,
\end{aligned}
\end{equation}
with $\alpha=\alpha_{12}e^{\alpha_{12}c_u\sqrt\varepsilon}$.
Furthermore, Assumption~\ref{ass_V} provides the following estimate:
$$
\begin{aligned}
\dot V&(x(t),u(t))= \frac{\partial V(x(t),u(t))}{\partial x}f(x(t),u(t))+\frac{\partial V(x(t),u(t))}{\partial u}\dot u(t)\\
&\le -\sigma_{21} V(x(t),u(t))+\sigma_{22}\sqrt{V(x(t),u(t))}\|\dot u(t)\|\\
&\le -\sigma_{21} V(x(t),u(t))+\frac{c_u\sigma_{22}}{\eta\sqrt\varepsilon}\sqrt{V(x(t),u(t))h(t,x(t))}.
\end{aligned}
$$
Given any $\sigma\in(0,\sigma_{21})$, let
$$
\eta_1(\varepsilon)=\frac{c_u L_h\sigma_{22}}{\sigma_{11}(\sigma_{21}-\sigma)\sqrt\varepsilon}.
$$
Then the last inequality together with~\eqref{est_h} implies that, for any $\varepsilon>0$, $\eta>\eta_1(\varepsilon)$,
$$
\begin{aligned}
\dot V&(x(t),u(t))\le - \sigma V(x(t),u(t))+\sqrt{V(x(t),u(t))}\varphi(t),
\end{aligned}
$$
where
$$
\begin{aligned}
\phi(t)=\frac{\alpha c_u\sigma_{22}}{\eta\sqrt\varepsilon}\Big(\|u^0&-u^*(0)\|+\nu\\
&+\frac{c_uL_h}{\sigma_{11}\eta\sqrt\varepsilon}\int\limits_0^t \sqrt{V(x(s),u(s))}ds\Big).
\end{aligned}$$
Solving the resulting comparison inequality yields the following estimate:
$$
\begin{aligned}
 &   \sqrt{V(x(t),u(t))}\le \sqrt{V(x^0,u^0)}e^{-\sigma t/2}\\
    &+\frac{1}{2}\int\limits_0^te^{-\sigma(t-s)/2}\phi(s)ds\\
    &\le\sqrt{V(x^0,u^0)}e^{-\sigma t/2}+\frac{\alpha c_u\sigma_{22}}{\sigma\eta\sqrt\varepsilon}\big(\|u^0-u^*(0)\|+2\nu \big)\\
    &+\frac{\alpha c_u^2 L_h\sigma_{22}}{2\sigma_{11}\eta^2\varepsilon}\int\limits_0^te^{-\sigma(t-s)/2}\int\limits_0^s \sqrt{V(x(p),u(p))}dpds.
\end{aligned}
$$
Applying integration by part to the last term, we continue the above estimate as
$$
\begin{aligned}
    \sqrt{V(x(t),u(t))}&\le \sqrt{V(x^0,u^0)}e^{-\sigma t/2}\\
    &+\frac{\alpha c_u\sigma_{22}}{\sigma\eta\sqrt\varepsilon}\big(\|u^0-u^*(0)\|+\nu\big)\\
    &+\frac{\alpha c_u^2 L_h\sigma_{22}}{\sigma\sigma_{11}\eta^2\varepsilon}\int\limits_0^t\sqrt{V(x(s),u(s))}ds.
\end{aligned}
$$
Using again the Gr\"onwall–Bellman inequality for $t\in[0,\eta\varepsilon]$ yields
$$
\begin{aligned}
    &\sqrt{V(x(t),u(t))}\le\frac{\alpha c_u\sigma_{22}}{\sigma\eta\sqrt\varepsilon}e^{\frac{\alpha c_u^2 L_h\sigma_{22}}{\sigma\sigma_{11}\eta}}\big(\|u^0-u^*(0)\|+\nu\big)\\
    &+\sqrt{V(x^0,u^0)}\Big(e^{-\sigma t/2}+\frac{\alpha c_u^2 L_h\sigma_{22}}{\sigma\sigma_{11}\eta^2\varepsilon}e^{\frac{\alpha c_u^2 L_h\sigma_{22}}{\sigma\sigma_{11}\eta}}\Big).
\end{aligned}
$$
Thus, for any $t\in[0,\eta\varepsilon]$,
\begin{equation}
    \label{estV1}
\begin{aligned}
   &\sqrt{V(x(t),u(t))}\le \frac{1}{\eta\sqrt\varepsilon}\Big({c_{v1}}\|u^0-u^*(0)\|+{\tilde\nu}\Big)\\
    &+\Big(e^{-\sigma t/2}+\frac{c_{v2}}{\eta^2\varepsilon}\Big)\sqrt{V(x^0,u^0)},
  \end{aligned}
\end{equation}
where $$c_{v1}=\frac{\alpha c_u\sigma_{22}}{\sigma}e^{\frac{\alpha c_u^2 L_h\sigma_{22}}{\sigma\sigma_{11}}\eta},\; c_{v2}=\frac{c_{v1}c_uL_h}{\sigma_{11}},\; \tilde \nu= c_{v1}\nu.$$
Consequently,
\begin{equation}
    \label{est_xu}
\begin{aligned}
   &\|x(t)-\ell(u(t))\|\le  \frac{1}{\alpha_{11}\eta\sqrt\varepsilon}\Big({c_{v1}}\|u^0-u^*(0)\|+{\tilde\nu}\Big)\\
    &+\frac{1}{\alpha_{11}}\Big(e^{-\sigma t/2}+\frac{c_{v2}}{\eta^2\varepsilon}\Big)\sqrt{V(x^0,u^0)},
  \end{aligned}
\end{equation}
for all $t\in[0,\eta\varepsilon]$.

\vskip1em

\emph{\textbf{Step 2.} Estimating the maximal deviations between the  $u$ component of the solution of system~\eqref{sys_dyn}--\eqref{cont_dyn} with $(x^0,u^0)\in\overline{B_{\delta_{x}}(\ell(u^*(t)))\times B_{\delta_u}(u^*(t))}$ and the solutions of the reduced system.}

As the next step, consider the auxiliary system
\begin{equation}\label{cont_dyn_au}
\dot {\bar u}_j=\frac{2\gamma}{\eta\sqrt\varepsilon}\sqrt{\pi  jh(t,\ell(\bar u))}\sin\left(\ln h(t_m,\ell(\bar u))+\frac{2\pi j t}{\eta\varepsilon}\right),
\end{equation}
with  $j=1,\dots,n_u$, $ t\in [0,\eta\varepsilon]$,   $\bar u(t^0)=u(t^0)\in D_u$.
Then
$$
\begin{aligned}
 \|&u(t)-\bar u(t)\|\le\frac{c_u L_{hx}}{\eta\sqrt\varepsilon}\int_{t_m}^{t_{m+1}}\|x(s)-\ell(\bar u(s))\|ds\\
 &\le\frac{c_u L_{hx}}{\eta\sqrt\varepsilon}\int_{t_m}^{t_{m+1}}\Big(\|x(s)-\ell( u(s))\|+L_\ell\|u(s)-\bar u(s)\|\Big)ds.\\
\end{aligned}
$$
Together with the Gr\"ownall--Bellman inequality, this implies
\begin{equation}
    \label{ubu}
\begin{aligned}
  \|u(t)-\bar u(t)\|\le \frac{c_{\tilde u 1}}{\eta}\|u^0&-u^*(0)\|\\
  &+\frac{c_{\tilde u 2}}{\eta\sqrt\varepsilon}\sqrt{V(x^0,u^0)}+\frac{c_{\tilde u 3}}{\eta}\tilde \nu
\end{aligned}
\end{equation}
for all $t\in[0,\eta\varepsilon]$,
where
$$c_{\tilde u 1}=c_{\nu1}c_{\tilde u 3} ,\; c_{\tilde u2}=\big(1+\dfrac{c_{\nu2}}{\eta}\Big)c_{\tilde u 3},\;     c_{\tilde u 3}=e^{c_u L_{hx} L_\ell \sqrt\varepsilon}\dfrac{c_u L_{hx}}{\alpha_{11}}.$$

\vskip1em

\emph{\textbf{Step 3.} Ensuring tracking properties for  the solutions of the reduced system.}

Based on the results of the paper~\cite{GE21}, we can derive the following representation of the solutions of system~\eqref{cont_dyn_au}:
$$
\bar u(\eta\varepsilon)=u^0-\varepsilon\gamma^2(u^0-u^*(0))+R_{\bar u},
$$
where $\|R_{\bar u}\|\le c_{R\bar u}\varepsilon^{3/2}$. Thus,
$$
\begin{aligned}
    \|\bar u(\eta\varepsilon)- u(\eta\varepsilon)\|\le (1&-\varepsilon\gamma^2)\|u^0-u^*(0)\|\\
    &+\|\bar u(\eta\varepsilon)-u^*(0)\|+c_{R\bar u}\varepsilon^{3/2}\\
    \le (1&-\varepsilon\gamma^2)\|u^0-u^*(0)\|+\nu+c_{R\bar u}\varepsilon^{3/2}.
\end{aligned}
$$
Given any $\rho>\nu$ and any $\varepsilon>0$, let us take any $\rho'\in(0,\rho)$, and  assume $\dfrac{\nu}{\rho'}<\varepsilon\gamma^2<1$.
Then, similarly to the approach of~\cite{GMZME18}, one can show that there exist an $\varepsilon_2>0$ and $\bar\gamma(\varepsilon)>0$ such that, for all $\varepsilon\in(0,\varepsilon_2)$, the solutions of system~\eqref{cont_dyn_au} with $\gamma=\bar\gamma(\varepsilon)$ satisfy the following properties:
\begin{itemize}
    \item if $\|u^0-u^*(0)\|\ge\rho'$, then $\|\bar u(\eta\varepsilon)-u^*(\eta\varepsilon)\|<\|u^0-u^*(0)\|$;
    \item if $\|u^0-u^*(0)\|<\rho'$, then $\|\bar u(\eta\varepsilon)-u^*(\eta\varepsilon)\|<\rho'$.
\end{itemize}
The rest of the proof is based on the estimate~\eqref{ubu} with a large  enough $\eta$ and follows the line of~\cite[Proof of Theorem~1, Step~4]{GE21}.

\end{document}